\documentclass{amsart}
\usepackage{amssymb}
\usepackage{amsfonts}
\usepackage{pstricks}
\usepackage{hyperref}

\theoremstyle{plain}

\renewcommand\bigskip\medskip

\def\N{\mathbb{N}}
\def\Z{\mathbb{Z}}

\def\Q{\mathbb{Q}}

\def\geq{\geqslant}

\begin{document}
\centerline{\bf{ON A FAMILY OF CONTINUED FRACTIONS IN $\Q((T^{-1}))$}}
\centerline{\bf{ASSOCIATED TO INFINITE BINARY WORDS}}
\centerline{\bf{ DERIVED FROM THE THUE-MORSE SEQUENCE}}

\vskip 0.5 cm
\centerline{\bf{by B. Allombert and A. Lasjaunias }}
\centerline {(Institut Math\'ematique de Bordeaux, France)}
\vskip 0.5 cm

\noindent Keywords Words: Combinatorics on words. Formal power series over $\Q$. Continued fractions.
\newline Mathematics Subject Classification 11J70, 68R15.
\vskip 0.5 cm
{\bf{Abstract.}} For each integer $n\geq 1$, we present an element  in $\Q((T^{-1}))$, having a power series expansion based on an infinite word $W(n)$, over the alphabet $\{+1,-1\}$ and whose continued fraction expansion has a particular pattern which is explicitly described. The word $W(1)$ is the Thue-Morse sequence and the following words are defined in a similar way.

\vskip 0.5 cm
\section{Introduction}
\par To any given infinite word $W=w_1w_2\dots w_n\dots$, whose letters take values in $\Q$, we associate an element in $\Q((T^{-1}))$ defined by $\sum_{n\geq 1}w_nT^{-n}$. In a previous work by the same authors [1], a particular infinite word over the alphabet $\{1,2\}$ was considered, leading as above to a generating function in $\Q((T^{-1}))$. The continued fraction expansion of this function was studied and it could be described partially. In a second and recent article [2], a full description of this continued fraction was given. Then, a question was arised about the existence of other infinite two-valued words leading to elements in $\Q((T^{-1}))$, with a particular and describable continued fraction. In this note, we answer this question positively by describing a family of such words and their corresponding continued fractions.
\newline For a general and basic information about continued fraction, particularly in power series fields, the reader may consult [3].   
\newline 
 The words considered here, belonging to a particular family,  have a simple pattern presented in the next section. In the last section we describe the continued fraction expansion in $\Q((T^{-1}))$, for each element $\theta$, attached as above to each $W$ in this family. The aim of this note is just to present these continued fractions, empirically discovered using a computer, consequently there will be no formal proofs. For this reason, to support this presentation, we invite the interested reader to consult a computer program ([4] and the link therein). 

\vskip 0.5 cm
\section{A family of infinite binary words} 
\par Let $W$ be a word over the alphabet $\{a,b\}$. We will use the following notation. If $W$ is a finite word and $k\geq 1$ is an integer we denote by $W^{[k]}=W,...,W$, the word obtained by concatenation of $k$ blocks equal to $W$. Moreover we let $W^*$ denote the word obtained by changing, in $W$, $a$ into $b$ and vice versa $b$ into $a$. 
\newline For $i\geq 1$, we define a sequence of finite words $(W(i)_n)_{n\geq 1}$ as follows :
$$W(i)_1=a \quad \text{and}\quad  W(i)_{n+1}=W(i)_n^{[i]},(W(i)_n^{[i]})^* \quad \text{for}\quad  n\geq 1.$$
Then we define the infinite word $W(i)$ as the projective limit of the sequence $(W(i)_n)_{n\geq 1}$, that is the word beginning by $W(i)_n$ for $n\geq 1$.
\newline Let us describe $W(1)$ and $W(2)$. The first one is the famous Thue-Morse sequence considered and studied independently by A. Thue and M. Morse about a century ago. Note that this same sequence was implicitely introduced, 175 years ago, by E. Prouhet.
\newline We have $W(1)_2=ab$ and $W(1)_3=abba$ and 
$W(1)=abbabaabbaababbaabb\dots$
\newline This sequence is registered in the on-line encyclopedia of integer sequences (OEIS, A010060). So is $W(2)$ the second one (A269723). 
\newline We have $$W(2)_2=aabb, \quad W(2)_3=aabbaabbbbaabbaa \quad \text{ and}$$ 
$$W(2)=aabbaabbbbaabbaaaabbaabbbbaabbaabbaa\dots$$ 
Note that, due to the defining relation, these words $W(i)$, for $i\geq 2$, can somehow be regarded as descending from the historical example $W(1)$. This will be underlined considering, in the next section, the continued fractions attached to them.     
 \vskip 0.5 cm
 \section{A family of infinite continued fractions in $\Q((T^{-1}))$} 

\par From now on, we choose the pair $(a,b)=(+1,-1)$. Hence, for $i\geq 1$, we consider the infinite word $W(i)=w_{i,1},w_{i,2},...,w_{i,n},...$ defined above where $w_{i,n}\in \{+1,-1\}$. Linked to this word, we define $\theta_i \in \Q((T^{-1}))$, by 
$$\theta_i=\sum_{n=1}^{\infty} w_{i,n}T^{-n}.$$ 
We shall describe the continued fraction expansion of $\theta_i$ for each $i\geq 1$. 
This description is made from several observations using a computer software (PARI/GP). These infinite continued fraction expansions are denoted as follows :
$$\theta_i=[0,a_{i,1},a_{i,2},\dots,a_{i,n},\dots]\quad \text{where}\quad a_{i,n}\in \Q[T].$$
In order to describe the partial quotients $a_{i,n}$, for each $i\geq 1$ and each $n\geq 1$, we will use two sequences. The first one, $(\lambda_{i,n})_{n\geq 1}$ in $\Q$, is the sequence of the leading coefficients of the partial quotients. The second one, $(b_{i,n})_{n\geq 1}$, is a sequence in $\Q[T]$, of unitary polynomials such that we have :
$$a_{i,n}=\lambda_{i,n}b_{i,n} \quad \text{for} \quad n\geq 1 \quad \text{ and for} \quad i\geq 1.$$
In our empirical search, we happened to discover first the case $i=2$. We could realise that the continued fraction was predictable, even though complicated enough. Then we tried the simpler case $i=1$ and we were surprised by the similarities with the second case. Later on, to understand better the general structure, we decided to look at the other cases. Several properties, concerning the first sequence $(\lambda_{i,n})_{n\geq 1}$ in $\Q$, are valid in all cases. Indeed, for all $i\geq 1$ and for all $k\geq 1$, we have observed the following three general equalities : 
$$(I):\quad \lambda_{i,4k-3}=\lambda_{i,2k-1},\qquad(II):\quad \lambda_{i,4k}=-\lambda_{i,4k-2}$$ 
$$\text{and}\quad (III):\quad \lambda_{i,4k-2}\lambda_{i,4k-1}\lambda_{i,4k}=\lambda_{i,2k}.$$
It is easy to obtain the first two values. For $i\geq 1$, we have $(\lambda_{i,1},\lambda_{i,2})=(1,1/2)$. Moreover, from $(I),(II)$ and $(III)$, we get easily $(\lambda_{i,3},\lambda_{i,4},\lambda_{i,5})=(-2,-1/2,-2)$. However, these three equalities do not allow to describe recursively the sequence $(\lambda_{i,n})_{n\geq 1}$ in $\Q$. To do so we will express $\lambda_{i,4k}$, using the following quantity :
$$u_{i,k}=\lambda_{i,4k}^{-1}-\lambda_{i,4k+4}^{-1} \quad \text{for }\quad i \geq 1  \quad \text{and }\quad k\geq 1.\eqno{(1)}$$
From $(1)$, since $\lambda_{i,4}^{-1}=-2$, we get easily
$$\lambda_{i,4k}=-(2+\sum_{j=1}^{j=k-1}u_{i,j})^{-1} \quad \text{for }\quad i \geq 1  \quad \text{and }\quad k\geq 2.\eqno{(2)}$$
We will see below that $u_{i,k}$ can be expressed in terms of preceding values of $\lambda_{i,j}$ and here will come a difference between the case $i=1$, attached to the historical word $W(1)$, and the other cases $i\geq 2$. Actually, the most striking difference, from the first case and the other ones, comes from the character bounded or not of the sequence of the degrees of the partial quotients for $\theta_i$.
\newline To be precise, we look first at the sequence $(b_{i,n})_{n\geq 1}$ in $\Q[T]$. This sequence, is very simple in the first case ($i=1$). Indeed, we observe that it is purely periodic of period $2$ and we have : 
$$ b_{1,n}=T-(-1)^n \qquad \text{for}\quad n\geq 1.\eqno{(IV)} $$
In the following cases ($i\geq 2$), this sequence of polynomials is more sophisticated and unbounded in degrees. In order to describe it, we need to introduce other sequences, for $i\geq 2$, of polynomials $(P_{i,m})_{m\geq 1}$ in $\Z[T]$. These sequences will be described below, but, with  this notation, for $i\geq 2$ and $k\geq 1$, we observed that :
$$b_{i,1}=T-1, \quad b_{i,2k}=(T^i-1)/(T-1) \quad \text{and}\quad b_{i,2k+1}=(T-1)P_{i,\nu_2(k)+1},\eqno{(V)}$$
 where $\nu_2(x)$ denotes the 2-adic valuation of an integer $x\in \N$. The sequence $(P_{i,m})_{m\geq 1}$ is defined recursively by the following relations :
 $$P_{i,1}=(T^{2i^2}-1)/(T^{2i}-1)\eqno{(\bf{P1})}$$ and
 $$ P_{i,m+1}(T)=P_{i,m}(T^{2i})+2(P_{i,m}(T^{2i})-P_{i,m}(1))/(T^{i}-1).\eqno{(\bf{P2})}$$
 Note that this definition of the polynomials $P_{i,m}$ implies the following formulas, concerning their degrees :
 $$\deg(P_{i,1})=2i(i-1) \quad \text{and }\quad \deg(P_{i,m+1})=2i\deg(P_{i,m}).$$
 \par We wish to make a brief digression about an arithmetical property of the elements $\theta_i$, for $i\geq 2$. Indeed, for $i\geq 2$ and $m\geq 1$, we have $\deg(P_{i,m})=(2i)^m(i-1)$. Hence, we have a precise knowledge of the sequence of the degrees of the partial quotients in the continued fraction expansion of $\theta_i$. From this sequence of degrees, one can compute the irrationality measure of the element. For more information on this matter the reader may consult [3, p. 11-13] or [1, p. 861-863]. In the present case, for $i\geq 2$, we obtain an irrationality measure for $\theta_i$ equal to $2i$. (Note that this is also true for $i=1$, since the partial quotients are bounded in degrees, in that case.) For $i\geq 2$, this proves that $\theta_i$ is a transcendental element in $\Q((T^{-1}))$, over $\Q(T)$.
 \vskip 0.5 cm
 \par We turn now back to the description of the sequences $(\lambda_{i,n})_{n\geq 1}$ of the leading coefficients of the partial quotients. 
 Here, again , there is a simplification in the first case $i=1$. In that case, we have observed that $u_{1,k}$ takes a simple form and we have
 $$u_{1,k}=2\lambda_{1,2k+1} \quad \text{for }\quad k \geq 1.\eqno{(3)}$$ 
 Combining $(2)$ and $(3)$, and recalling that $\lambda_{1,1}=1$, we get
 $$\lambda_{i,4k}=-(2\sum_{j=0}^{j=k-1}\lambda_{1,2j+1})^{-1} \quad \text{for }\quad k \geq 1 .\eqno{(4)}$$
  \par Combining the formulas $(I),(II),(III),(IV)$ and $(4)$, we can describe fully the continued fraction expansion for $\theta_1$, corresponding to the historical case of the word $W(1)$. So we have the following unproved statement (but empiricaly verified):
  \vskip 0.5 cm
 \par \emph{ Let $(\lambda_n)_{n\geq 1}$ be the sequence of rational numbers defined recursively by $\lambda_1=1$, $\lambda_2=1/2$ and for $k\geq 1$ : 
  $\lambda_{4k-3}=\lambda_{2k-1}$, $\lambda_{4k}=-\lambda_{4k-2}$, 
   $$\lambda_{4k-2}= (2\sum_{j=0}^{j=k-1}\lambda_{2j+1})^{-1}\quad \text{ and }\quad \lambda_{4k-1}=-\lambda_{4k-2}^{-2}\lambda_{2k}.$$
   Then, in $\Q((T^{-1}))$, we have the following identity
   $$[0,T+1,(1/2)(T-1),\dots,\lambda_n(T-(-1)^n),\dots ]=\sum_{i\geq 1}w_iT^{-i},$$
   where $(w_i)_{i\geq 1}$ is the Thue-Morse sequence over the alphabet $\{+1,-1\}$.}
   \vskip 0.5 cm
  \par Concerning the cases $i\geq 2$, we saw above the description of the sequence $(b_{i,n})_{n\geq 1}$ in $\Z[T]$. Contrarily  to the first case, this sequence $(b_{i,n})_{n\geq 1}$ is interfering into the sequence $(\lambda_{i,n})_{n\geq 1}$. Indeed, this happens through the quantities $u_{i,j}$, $(3)$ being replaced by :
  $$u_{i,j}=2\lambda_{i,2k+1}P_{i,\nu_2(j)+1}(1) \quad \text{for }\quad j \geq 1 \quad \text{and }\quad i\geq 2.\eqno{(5)}$$ 
  Consequently, joining the three equalities $(I),(II),(III)$ and also $(2),(5)$ the sequence $(\lambda_{i,n})_{n\geq 1}$ is defined recursively. Hence, together with $(V)$, we have a full description of the continued fraction for $\theta_i$, where $i\geq 2$.
  \newline At last, we make a brief remark concerning the values $P_{i,m}(1)$, appearing in $(5)$. These values are obviously coming from the recurrent definition of the sequence $(P_{i,m})_{m\geq 1}$ in $\Z[T]$. Indeed, from {(\bf{P1})} and {(\bf{P2})}, we get
  $$P_{i,1}(1)=i \quad \text{and }\quad P_{i,m+1}(1)=P_{i,m}(1) + 4P'_{i,m}(1) \quad \text{for }\quad m \geq 1.$$
  Basic computations show that $P'_{i,1}(1)=i^2(i-1)$ and $P'_{i,2}(1)=(8/3)i^4(i-1)(2i-1)$. While for a general value of the pair $(i,m)$, a closed-form expression for $P'_{i,m}(1)$ seems difficult to write. However these values can be obtained, by help of a computer, for all $i\geq 1$ and $m\geq 1$. For this, the reader may look at the computer program indicated at the end of the introduction [4].

Bill.Allombert@math.u-bordeaux.fr   \hskip 1 cm  lasjauniasalain@gmail.com   
\end{document}